\documentclass[preprint,12pt, a4paper]{elsarticle}

\usepackage{hyperref}

\journal{SoftwareX}

\usepackage{color}
\usepackage{amsmath}
\usepackage{graphicx}
\usepackage{subfigure}
\usepackage{rotating}
\usepackage{multirow, makecell}
\usepackage{geometry}
\usepackage{appendix}
\usepackage{float}

\biboptions{sort&compress}              

\begin{document}

\sloppy
	
\begin{frontmatter}
		
\title{MSAT: Matrix stability analysis tool for shock-capturing schemes}

\author[mymainaddress]{Weijie Ren}

\author[mymainaddress]{Wenjia Xie\corref{mycorrespondingauthor}}
\cortext[mycorrespondingauthor]{Corresponding author}
\ead{xiewenjia@nudt.edu.cn}

\author[mymainaddress]{Ye Zhang}

\author[mymainaddress]{Hang Yu}

\author[mymainaddress]{Zhengyu Tian}

\address[mymainaddress]{College of Aerospace Science and Engineering, National University of Defense Technology, Hunan 410073, China}

\begin{abstract}
	The simulation of supersonic or hypersonic flows often suffers from numerical shock instabilities if the flow field contains strong shocks, limiting the further application of shock-capturing schemes. In this paper, we develop the unified matrix stability analysis method for schemes with three-point stencils and present MSAT, an open-source tool to quantitatively analyze the shock instability problem. Based on the finite-volume approach on the structured grid, MSAT can be employed to investigate the mechanism of the shock instability problem, evaluate the robustness of numerical schemes, and then help to develop robust schemes. Also, MSAT has the ability to analyze the practical simulation of supersonic or hypersonic flows, evaluate whether it will suffer from shock instabilities, and then assist in selecting appropriate numerical schemes accordingly. As a result, MSAT is a helpful tool that can investigate the shock instability problem and help to cure it.
\end{abstract}

\begin{keyword}
	Shock instability, Matrix stability analysis, MUSCL, ROUND
\end{keyword}
		
\end{frontmatter}

\begin{table}[H]
	\centering
	\footnotesize 
	\begin{tabular}{|l|p{6.2cm}|p{7.8cm}|}
	\hline
	\textbf{Nr.} & \textbf{Code metadata description} & \textbf{Please fill in this column} \\
	\hline
	C1 & Current code version & v1.0.0 \\
	\hline
	C2 & Permanent link to code/repository used for this code version &https://github.com/JameRwj/MSAT.git\\
	\hline
	C3  & Permanent link to Reproducible Capsule & \\
	\hline
	C4 & Legal Code License   & GNU General Public License v3.0 \\
	\hline
	C5 & Code versioning system used & git \\
	\hline
	C6 & Software code languages, tools, and services used &Fortran90, MKL, Lapack\\
	\hline
	C7 & Compilation requirements, operating environments \& dependencies & \scriptsize{Test on Windows with Visual Studio 2022 and Intel oneAPI 2023; Test on Linux with gfortran 9.4.0 and Lapack 3.8.0.} \\
	\hline
	C8 & If available Link to developer documentation/manual &  \\
	\hline
	C9 & Support email for questions &xiewenjia@nudt.edu.cn \\
	\hline
	\end{tabular}
	\caption{Code metadata}
\end{table}

\section{Motivation and significance}\label{section 1}
With decades of development, computational fluid dynamics (CFD) has been widely used in scholarly research and industrial development. Unfortunately, when performing the simulation of supersonic and hypersonic flows, the flow field computed by the modern shock-capturing schemes is often characterized by numerical shock instabilities. As one of the most famous kinds of shock anomalies, the carbuncle phenomenon was first observed by Perry and Imlay \cite{Perry1988} when they simulated the supersonic flow around the blunt-body with Roe Riemann solver \cite{Roe1981}. The carbuncle phenomenon refers to the spurious solution of blunt-body calculations in which a protuberance grows ahead of the bow shock along the stagnation line \cite{Quirk1994}. Also, it has been generally considered to be a concrete manifestation of the numerical shock instability problem.

When the shock instability phenomenon happens, the flow field will exhibit obvious unphysical characteristics, resulting in simulation results that are deemed unreliable. Consequently, it is essential to investigate the primary numerical characteristics and the underlying mechanism of the numerical shock instability phenomenon and then explore possible methods to cure the shock instability. To this end, tools for analyzing the shock instability problem are needed. Quirk \cite{Quirk1994} first conducts a comprehensive study on the shock instability and proposes the method of linearized perturbation analysis. The linearized perturbation analysis is based on the odd-even decoupling, and the main idea is to examine how different Riemann solvers evolve the sawtooth-type initial data. Quirk's linearized algorithm provides a useful linearized temporal evolution model for the perturbations of different variables, as well as their mutual interactions, thus making it a useful tool for assessing the shock stability of various Riemann solvers. As a result, Quirk's linearized algorithm and its alterations are widely employed to investigate the internal mechanism of the shock instability problem \cite{Quirk1994,Pandolfi2001,Gressier2005,Xie2017,Xie2021}. However, the linearized perturbation analysis, performed manually instead of through programming, can only provide qualitative analysis. Moreover, it has been demonstrated that shock instability is affected by a variety of issues, such as Riemann solvers, shock intensity, computational grid, numerical shock structure, and so on \cite{Quirk1994,Kitamura2009,Henderson2007,Xie2017,Pandolfi2001,Gressier2005}. However, Quirk's linearized algorithm deals only with the perturbations superimposed on a steady mean flow without shocks and does not consider the effects of other issues. Except for the linearized perturbation analysis, the numerical experiment is also an effective way to study the shock instability problem. By performing various well-designed numerical experiments, the shock stability with different conditions can be investigated. Readers are referred to references \cite{Tu2014,Ohwada2013,Kitamura2009,Kitamura2010,Kitamura2013a,Xie2017,Xie2021} for detailed applications of the numerical experiment. Although the numerical experiment is easy to implement and can involve the impacts of different factors, it is also unable to carry out the quantitative analysis.

To investigate the shock instability problem quantitatively, Dumbser et al.\cite{Dumbser2004} propose the matrix stability analysis method, coupling the capacity of schemes to stably capture shocks with the eigenvalues of the stability matrix. Based on this method, the evolution of the perturbation errors can be quantitatively analyzed easily by programming. In addition, it allows incorporating the features of the numerical shock, effects of the computational grid, and boundary conditions \cite{Simon2018a}. As a result, the matrix stability analysis is widely employed to investigate the mechanism of the shock instability\cite{Dumbser2004,Chauvat2005,Shen2014,Liu2020,Chen2018a} and validate the shock stability of the novel shock-stable schemes\cite{Chen2018,Xie2019,Hu2022,Chen2018c,Chen2023,Sun2022a}. 

Unfortunately, no open-source resources are currently available to easily implement the matrix stability analysis method. As a result, if researchers want to utilize this method to investigate the shock instability problem, they must develop the framework of the matrix analysis from scratch, which is time-consuming and error-prone. Furthermore, the matrix stability analysis method proposed by Dumbser et al.\cite{Dumbser2004} is only applicable to the first-order scheme. Since higher-order schemes are applied more frequently in practical simulations of supersonic or hypersonic flows, it is necessary to investigate the shock instability problem for them. Results in \cite{Pandolfi2001,Liou2000,Zhang2017b} show that the stability of the second-order scheme may differ from the first-order case, and the limiter used in the second-order scheme plays a vital role in the shock instability. It also has been demonstrated that fifth-order schemes are at higher risk of the shock instability \cite{Tu2014,Jiang_Effective_2017}. However, there have been few efforts to analyze the shock instability problem for high-order schemes. Such a situation is mainly due to the lack of an effective analytical tool. In the previous work \cite{2305.03281}, we develop the matrix stability analysis method for the finite-volume MUSCL approach. Moreover, in the current study, based on the work in \cite{2305.03281}, we develop a unified matrix stability analysis method for schemes with three-point stencils and present an open-source tool, MSAT, for analyzing the shock instability of these schemes, such as the second-order MUSCL scheme and ROUND scheme \cite{VanLeer1979,VanLeer2021}.

The outline of the rest of this paper is as follows. Section \ref{section 2} describes the matrix stability analysis method for schemes with three-point stencils and the architecture of MSAT. The quantitative validation of MSAT and several illustrative examples are presented in section \ref{section 3}. Furthermore, the impact of MSAT is discussed in section \ref{section 4}.

\section{Software description}\label{section 2}
\subsection{Governing equations and finite-volume discretization}\label{subsection 2.1}
In the current study, we employ the two-dimensional Euler equations, which can be written in the integral form as
\begin{equation}\label{eq Euler equations}
  \frac{\partial}{\partial t} \int_{\Omega} \mathbf{U} {\rm{d}}{\Omega} + \oint_{\partial \Omega} {\mathbf{F}} {\rm{d}}S = 0,
\end{equation}
where $ \mathbf{U} $ denotes the vector of conservative variables and ${\mathbf{F}}$ is the flux component normal to $\partial \Omega$, which is the boundary of the control volume $\Omega$. We consider discretizing the system (\ref{eq Euler equations}) with the cell-centered finite-volume method over the 2D domain subdivided into structured quadrilateral cells, which can be written in the following forms 
\begin{equation}\label{eq discrete Euler equations}
  \frac{\mathrm{d} {\mathbf{U}}_{i,j}}{\mathrm{~d} t}=-\frac{1}{\left|\Omega_{i,j}\right|} \sum_{k=1}^{4} \mathcal{L}_{k} \mathbf{F}_{k},
\end{equation}
where $ {\mathbf{U}}_{i,j} $ denotes the average of \textbf{U} on $ \Omega _{i,j} $. $ \left|\Omega_{i,j}\right| $ is the volume of $ \Omega _{i,j} $ and $\mathcal{L}_{k}$ stands for the length of the cell interface. $\mathbf{F}_{k}$ is the numerical flux and is the function of the variables on the left and right sides of the interface. When considering the interface between $ \Omega_{i,j} $ and $ \Omega_{i+1,j} $, the numerical flux can be written as
\begin{equation}\label{eq flux function}
  \mathbf{F}_{i+1 / 2, j}=\mathbf{F}_{i+1 / 2, j}\left(\mathbf{U}_{i+1 / 2, j}^L, \mathbf{U}_{i+1 / 2, j}^R\right).
\end{equation}
where $ \mathbf{U}_{i+1 / 2,j}^{L} $ and $ \mathbf{U}_{i+1 / 2,j}^{R} $ are the variables on the left and right sides of the interface. If the numerical scheme is first-order accurate, only one point is used to determine the values of $\mathbf{U}_{i+1/2,j}^{L}$ and $\mathbf{U}_{i+1/2,j}^{R}$, which can be expressed as follows
\begin{equation}  
	\mathbf{U}_{i+1 / 2,j}^{L} =  \mathbf{U}_{i,j}\enspace \text{and} \enspace \mathbf{U}_{i+1 / 2,j}^{R} = \mathbf{U}_{i+1,j}.
\end{equation}
To achieve higher-order accuracy, more points are employed to reconstruct $\mathbf{U}_{i+1/2,j}^{L}$ and $\mathbf{U}_{i+1/2,j}^{R}$. In this peper, we concrete on the schemes with three-point stencils, such as the finite-volume MUSCL approach \cite{VanLeer1979,VanLeer2021} and ROUND scheme \cite{Deng_Unified_2023,Cheng_Accurate_2023,Deng_Largeeddy_2023}. According to the MUSCL approach, $ \mathbf{U}_{i+1 / 2,j}^{L} $ and $ \mathbf{U}_{i+1 / 2,j}^{R} $ can be written as
\begin{equation}\label{eq MUSCL method}
	\begin{aligned}
		& \mathbf{U}_{i+1 / 2,j}^{L}=\mathbf{U}_{ i,j}+\frac{1}{2} \Psi_{i+1 / 2,j}^{L}\left(\mathbf{U}_{i,j}-\mathbf{U}_{i-1,j}\right) \\
		& \mathbf{U}_{i+1 / 2,j}^{R}=\mathbf{U}_{i+1,j}-\frac{1}{2} \Psi_{i+1 / 2,j}^{R}\left(\mathbf{U}_{i+2,j}-\mathbf{U}_{i+1,j}\right)
	\end{aligned},
\end{equation}
where $ \Psi_{i+1/2,j}^{L/R} $ is the limiter function. The form of the ROUND scheme used in this paper is
\begin{equation}  \label{eq ROUND}
	\hat{\mathbf{U}}_{i+1/2,j}^L=\begin{cases}
		\min\left\{\left(\frac{1}{3}+\frac{5}{6}\hat{\mathbf{U}}_{i,j}^{L}\right)\omega_{0}+2.0\hat{\mathbf{U}}_{i,j}^{L}\left(1-\omega_{0}\right),2.0\hat{\mathbf{U}}_{i,j}^{L}\right\}&0.0<\hat{\mathbf{U}}_{i,j}^{L}\leq0.5,\\
		\begin{aligned}
			\min\left\{\left(\frac{1}{3}+\frac{5}{6}\hat{\mathbf{U}}_{i,j}^{L}\right)\omega_{1}+\left(\lambda_{1}\hat{\mathbf{U}}_{i,j}^{L}-\lambda_{1}+1.0\right)\left(1-\omega_{1}\right)\right.,\\
			\left.\lambda_{1}\hat{\mathbf{U}}_{i,j}^{L}-\lambda_{1}+1.0\right\}\
		\end{aligned}&0.5<\hat{\mathbf{U}}_{i,j}^{L}\leq1.0,\\
		\hat{\mathbf{U}}_{i,j}^{L}&\text{otherwise},
	\end{cases}
\end{equation}
where
\begin{equation}  \label{eq Normalised Reconstruction Value}
	\hat{\mathbf{U}}_{i+1/2,j}^L=\frac{\mathbf{U}_{i+1/2,j}^L-\mathbf{U}_{i-1,j}}{\mathbf{U}_{i+1,j}-\mathbf{U}_{i-1,j}}
	\qquad
	\hat{\mathbf{U}}_{i,j}^L=\frac{\mathbf{U}_{i,j}-\mathbf{U}_{i-1,j}}{\mathbf{U}_{i+1,j}-\mathbf{U}_{i-1,j}}.
\end{equation}
It can be found that both the MUSCL approach and ROUND scheme use three-point stencils to reconstruct $\mathbf{U}_{i+1/2,j}^{L}$ and $\mathbf{U}_{i+1/2,j}^{R}$ and they can be expressed in a unified form
\begin{equation}   \label{eq unified expression}
	\begin{aligned}
		\mathbf{U}_{i+1/2,j}^L&=\boldsymbol{\alpha}_{i-1,j}^L\mathbf{U}_{i-1,j}+\boldsymbol{\alpha}_{i,j}^L\mathbf{U}_{i,j}+\boldsymbol{\alpha}_{i+1,j}^L\mathbf{U}_{i+1,j},\\
		\mathbf{U}_{i+1/2,j}^R&=\boldsymbol{\alpha}_{i,j}^R\mathbf{U}_{i,j}+\boldsymbol{\alpha}_{i+1,j}^R\mathbf{U}_{i+1,j}+\boldsymbol{\alpha}_{i+2,j}^R\mathbf{U}_{i+2,j}.
	\end{aligned}
\end{equation}

\subsection{The matrix stability analysis method for schemes with three-point stencils}\label{subsection 2.2}
In our previous work \cite{2305.03281}, the matrix stability analysis method for the second-order MUSCL scheme is proposed. And in this section we further extend this method to the schemes with three-point stencils. For the stability analysis of a steady field, it can be assumed as follows
\begin{equation}\label{eq variables decomposition}
	\mathbf{U}_{i,j}=\mathbf{U}_{i,j}^{0}+\delta \mathbf{U}_{i,j},
\end{equation}
where $ \mathbf{U}_{i,j}^{0} $ denotes the steady mean value and $ \delta \mathbf{U}_{i,j} $ is the small numerical random perturbation. Substituting (\ref{eq variables decomposition}) into (\ref{eq unified expression}), we can get
\begin{equation}  
	\begin{aligned}
		\delta\mathbf{U}_{i+1/2,j}^L&=\boldsymbol{\alpha}_{i-1}^L\delta\mathbf{U}_{i-1,j}+\boldsymbol{\alpha}_{i}^L\delta\mathbf{U}_{i,j}+\boldsymbol{\alpha}_{i+1}^L\delta\mathbf{U}_{i+1,j}\\
		\delta\mathbf{U}_{i+1/2,j}^R&=\boldsymbol{\alpha}_{i}^R\delta\mathbf{U}_{i,j}+\boldsymbol{\alpha}_{i+1}^R\delta\mathbf{U}_{i+1,j}+\boldsymbol{\alpha}_{i+2}^R\delta\mathbf{U}_{i+2,j}
	\end{aligned} 
\end{equation}
According to (\ref{eq flux function}), the numerical flux can be linearized around the steady mean value as
\begin{equation}\label{eq simplify second-order flux linearized}
	\begin{aligned}
		\mathbf{F}_{i+1/2,j}& = \mathbf{F}_{i+1/2,j}\left( \mathbf{U}^{L,0}_{i+1/2,j},\mathbf{U}^{R,0}_{i+1/2,j}\right)\\
		&+\boldsymbol{\beta}_{i+1/2,j}^{i+2,j}\delta\mathbf{U}_{i+2,j}+\boldsymbol{\eta}_{i+1/2,j}^{i+1,j}\delta\mathbf{U}_{i+1,j}+\boldsymbol{\eta}_{i+1/2,j}^{i,j}\delta\mathbf{U}_{i,j}+\boldsymbol{\beta}_{i+1/2,j}^{i-1,j}\delta\mathbf{U}_{i-1,j}
	\end{aligned},
\end{equation}
where 
\begin{equation}  
	\begin{aligned}
		\boldsymbol{\beta}_{i+1/2,j}^{i+2,j} &= \frac{\partial\mathbf{F}_{i+1/2,j}}{\partial\mathbf{U}_{i+1/2,j}^R}\boldsymbol{\alpha}_{i+2,j}^R\\
		\boldsymbol{\eta}_{i+1/2,j}^{i+1,j} &= \frac{\partial\mathbf{F}_{i+1/2,j}}{\partial\mathbf{U}_{i+1/2,j}^L}\boldsymbol{\alpha}_{i+1,j}^L+\frac{\partial\mathbf{F}_{i+1/2,j}}{\partial\mathbf{U}_{i+1/2,j}^R}\boldsymbol{\alpha}_{i+1,j}^R\\
		\boldsymbol{\eta}_{i+1/2,j}^{i,j} &= \frac{\partial\mathbf{F}_{i+1/2,j}}{\partial\mathbf{U}_{i+1/2,j}^L}\boldsymbol{\alpha}_{i,j}^L+\frac{\partial\mathbf{F}_{i+1/2,j}}{\partial\mathbf{U}_{i+1/2,j}^R}\boldsymbol{\alpha}_{i,j}^R\\
		\boldsymbol{\beta}_{i+1/2,j}^{i-1,j} &= \frac{\partial\mathbf{F}_{i+1/2,j}}{\partial\mathbf{U}_{i+1/2,j}^L}\boldsymbol{\alpha}_{i-1,j}^L
	\end{aligned}.
\end{equation}
Substituting (\ref{eq simplify second-order flux linearized}) and (\ref{eq variables decomposition}) into (\ref{eq discrete Euler equations}), the linear error evolution mode can be obtained as follows
\begin{equation}
	\begin{aligned}\label{eq second-order linear error evolution model for conservative variables}
		\frac{\mathrm{d} \delta \mathbf{U}_{i, j}}{\mathrm{dt}}=&-\left(\boldsymbol{\mu}_{i+1/2,j}^{i,j}+\boldsymbol{\mu}_{i,j+1/2}^{i,j}+\boldsymbol{\mu}_{i-1/2,j}^{i,j}+\boldsymbol{\mu}_{i,j-1/2}^{i,j}\right)\delta\mathbf{U}_{i,j}\\
		&-\left(\boldsymbol{\mu}_{i+1/2,j}^{i+1,j}+\boldsymbol{\xi}_{i-1/2,j}^{i+1,j}\right)\delta\mathbf{U}_{i+1,j}-\left(\boldsymbol{\mu}_{i,j+1/2}^{i,j+1}+\boldsymbol{\xi}_{i,j-1/2}^{i,j+1}\right)\delta\mathbf{U}_{i,j+1}\\
		&-\left(\boldsymbol{\mu}_{i-1/2,j}^{i-1,j}+\boldsymbol{\xi}_{i+1/2,j}^{i-1,j}\right)\delta\mathbf{U}_{i-1,j}-\left(\boldsymbol{\mu}_{i,j-1/2}^{i,j-1}+\boldsymbol{\xi}_{i,j+1/2}^{i,j-1}\right)\delta\mathbf{U}_{i-1,j}\\
		&-\boldsymbol{\xi}_{i+1/2,j}^{i+2,j}\delta\mathbf{U}_{i+2,j}-\boldsymbol{\xi}_{i,j+1/2}^{i,j+2}\mathbf{U}_{i,j+2}-\boldsymbol{\xi}_{i-1/2,j}^{i-2,j}\delta\mathbf{U}_{i-2,j}-\boldsymbol{\xi}_{i,j-1/2}^{i,j-2}\delta\mathbf{U}_{i,j-2}
	\end{aligned},
\end{equation}
where
\begin{equation}
	\xi_{i\pm1/2, j\pm1/2}=\frac{\mathcal{L}_{i\pm1/2, j\pm1/2}}{|\Omega_{i,j}|}\beta_{i\pm1/2, j\pm1/2}\quad , \quad \mu _{i\pm1/2, j\pm1/2}=\frac{\mathcal{L}_{i\pm1/2, j\pm1/2}}{|\Omega _{i,j}|}\eta _{i\pm1/2, j\pm1/2}.
\end{equation}
It can be found that (\ref{eq second-order linear error evolution model for conservative variables}) holds for all cells in the computational domain and we can get \cite{Dumbser2004}
\begin{equation}\label{eq linear error evolution model for all cells}
	\frac{\mathrm{d}}{\mathrm{d} t}\left(\begin{array}{c}
		\delta \mathbf{U}_{1,1} \\
		\vdots \\
		\delta \mathbf{U}_{imax,jmax}
		\end{array}\right)=\mathbf{S} \cdot\left(\begin{array}{c}
		\delta \mathbf{U}_{1,1} \\
		\vdots \\
		\delta \mathbf{U}_{imax,jmax}
		\end{array}\right),
\end{equation}
where \textbf{S} is called the stability matrix in the present study. When considering the evolution of the initial error only, the solution of (\ref{eq linear error evolution model for all cells}) is
\begin{equation}\label{eq solution of linear error evolution model}
	\left(\begin{array}{c}
		\delta \mathbf{U}_{1,1} \\
		\vdots \\
		\delta \mathbf{U}_{imax,jmax}
		\end{array}\right)(t)=\mathrm{e}^{\mathbf{S} t} \cdot\left(\begin{array}{c}
		\delta \mathbf{U}_{1,1} \\
		\vdots \\
		\delta \mathbf{U}_{imax,jmax}
		\end{array}\right)_{t=0}.
\end{equation}
It can be demonstrated that the solution (\ref{eq solution of linear error evolution model}) will be bounded if all the eigenvalues of \textbf{S} is negative \cite{Dumbser2004}, as a result of which we can get the stability criterion
\begin{equation}\label{eq stability criterion}
	\max (\operatorname{Re}(\lambda(\mathbf{S}))) \leq 0.
\end{equation}

Note that equation (\ref{eq simplify second-order flux linearized}) is accurate only when the numerical flux is differentiable at the mean value, which is not always holding, for example, when the Roe solver is employed and the shock is exactly between two cells ($ \varepsilon =0 \enspace \text{or} \enspace 1 $). So, if Roe solver is analyzed, numerical shock structure must exist ($ 0<\varepsilon<1 $). The gradients of the numerical flux functions {such as $ \dfrac{\partial \mathbf{F}_{i+1/2,j}}{\partial\mathbf{U}_{i+1/2,j}^{L,0}} $} can be calculated as follows
\begin{equation}
	\dfrac{\partial \mathbf{F}_{i+1/2,j}}{\left(\partial \mathbf{U}_{i+1/2,j}^{L,0}\right)_k}=\dfrac{\mathbf{F}_{i+1/2,j}(\mathbf{U}_{i+1/2,j}^{L,0}+\delta\mathbf{I}_k,\mathbf{U}_{i+1/2,j}^{R,0})-\mathbf{F}_{i+1/2,j}(\mathbf{U}_{i+1/2,j}^{L,0}-\delta\mathbf{I}_k,\mathbf{U}_{i+1/2,j}^{R,0})}{2\delta},
\end{equation}
where $ \mathbf{I}_k $ is unit vector of which the $ kth $ component is 1, and $ \delta =10^{-7} $ \cite{Dumbser2004,Shen2014}.

By the method mentioned in this section, the shock instability of the scheme with three-point stencils is related to $ \max (\operatorname{Re}(\lambda(\mathbf{S}))) $, and can be analyzed quantitatively. It should be noted that since (\ref{eq solution of linear error evolution model}) is the analytical solution of (\ref{eq linear error evolution model for all cells}), there is no need to use temporal discretization schemes in this method.

\subsection{Software architecture}
One of the main purposes of MSAT is to study the mechanism of carbuncle phenomenon or shock instability. The carbuncle phenomenon is conventionally referred to as the distorted shock ahead of the blunt-body in the supersonic or hypersonic flow \cite{Quirk1994} and usually occurs in the flow field with strong shocks in practical flow simulations. However, it is necessary to choose a sample and representative test case when investigating its mechanism. MSAT is based on the test case of the 2D steady normal shock since it is easy to analyze and shares the fundamental characteristics of the carbuncle phenomenon. It has been demonstrated that if a scheme produces unstable solutions for the 2D steady normal shock problem, then it will also suffer from the blunt-body carbuncle \cite{Ismail2006,Dumbser2004,Kitamura2009}. As a result, the 2D steady normal shock problem has been widely used in studying the shock instability problem and developing the robust scheme \cite{Dumbser2004,Simon2018a,Kitamura2009,Kitamura2010,Kitamura2013a,Xie2017}. Meanwhile, MSAT can also be used to analyze whether the practical simulation of supersonic/hypersonic flows will suffer from the shock instability problem, if the initial flow is given. All codes in MSAT are written in Fortran, and the architecture is shown in Fig.\ref{fig Software architecture diagram}. As shown, a brief description is given as follows:

\begin{figure}[htbp]
	\centering
	\includegraphics[width=0.9\textwidth]{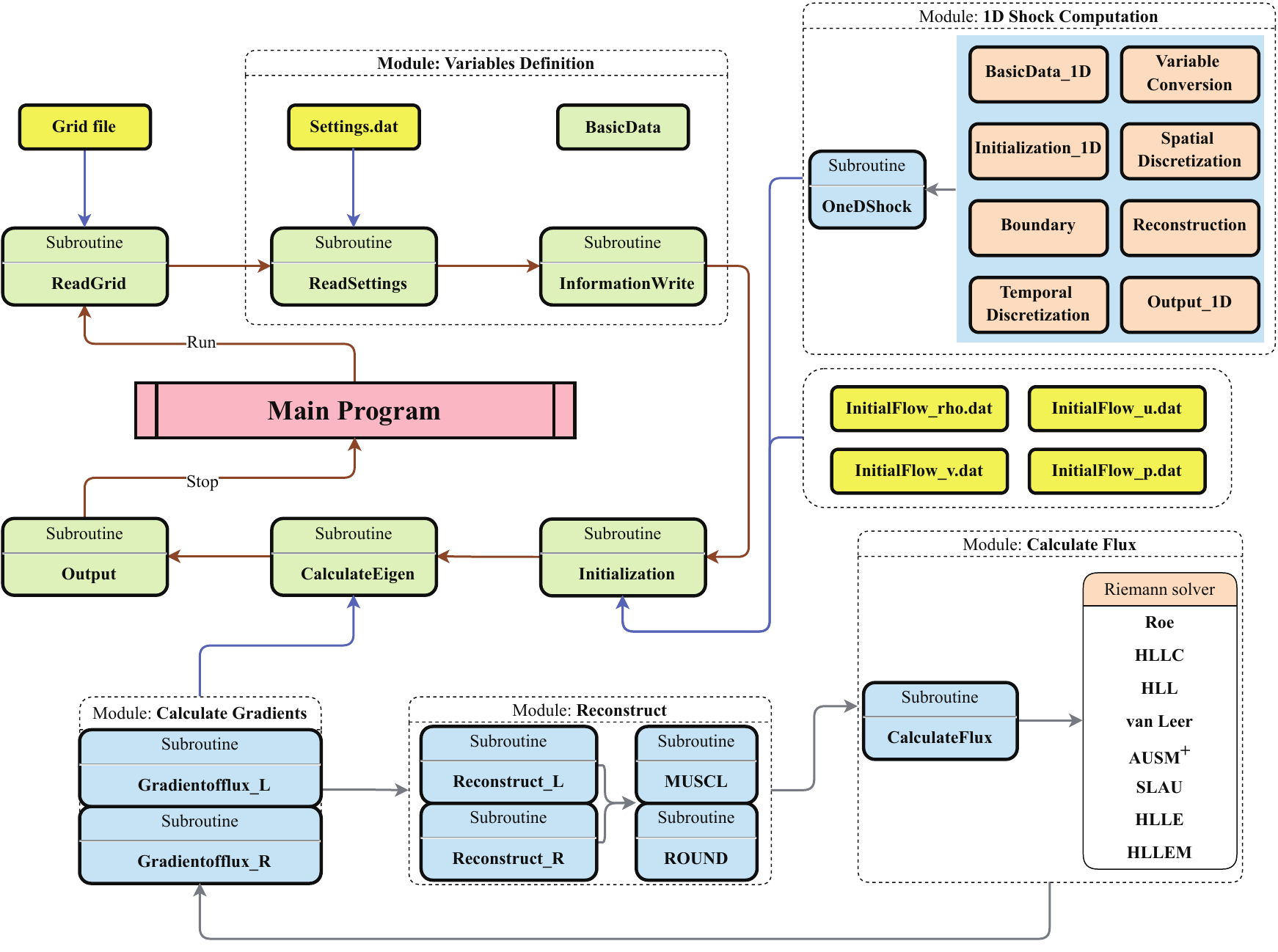}
	\caption{Software architecture diagram.}
	\label{fig Software architecture diagram}
\end{figure}

\begin{itemize}
  \item ``\textit{Main Program}", as the name suggests, is the main program of MSAT. It calls other subroutines to run the software.
  
  \item The role of subroutine ``\textit{ReadGrid}" is to read the coordinates of the grid nodes stored in the grid file. The computational grid should be structured, with the origin of the coordinates in the lower-left corner. There are two parts in the grid file. The first part specifies the number of grid nodes in the $ x $ and $ y $ dimensions. The second part is divided into three columns: the first column contains the $ x $ coordinates of each grid point, the second column contains the $ y $ coordinates, and the third column displays the $ z $ values, which should be 0 due to the two-dimensional nature of the computation.
  
  \item The module ``\textit{Variables Definition}'' defines, reads, and shows the main settings of the computation, which contains three parts: ``\textit{BasicData}'', ``\textit{ReadSettings}'', and ``\textit{InformationWrite}''. Thereinto, the global variables are defined in ``\textit{BasicData}'' and the main settings of the analysis are read by ``\textit{ReadSettings}'' from ``\textit{Settings.dat}'', including the reconstruction method, limiter function, Riemann solver, test case, numerical shock structure, Mach number, the way to initialize the flow field, and the iteration steps of 1D computation. Finally, the settings and grid information will be displayed on the screen by the subroutine ``\textit{InformationWrite}''.
  
  \item As mentioned in section \ref{subsection 2.2}, the matrix stability analysis is based on the steady mean value. So before calculation, there must be a stable flow field, which is performed by the subroutine ``\textit{Initialization}''. The initialization method of the 2D steady normal shock problem differs from that of the other cases.
  
  \begin{itemize} 
	\item For the 2D steady normal shock problem, the initial flow field is obtained by MSAT itself. There are two methods to initialize the flow field. The first way is depending on the Rankine-Hugoniot conditions across the normal shock. Readers can read \cite{Xie2021,2305.03281} for detailed information. Also, the 2D flow field can be initialized by projecting the steady flow field from 1D computation onto the 2D domain, which is also employed in \cite{Dumbser2004}. The 1D steady flow field is computed by the module ``\textit{1D shock computation}''.
	\item For other test cases, the initial flow field needs to be computed by other programs and is stored in ``\textit{InitialFlow\_rho.dat}'', ``\textit{InitialFlow\_u.dat}'', ``\textit{InitialFlow\_v.dat}'', ``\textit{InitialFlow\_p.dat}'', respectively. It should be noted that the variables of the flow field should be stored in the cell center and consistent with the distribution of the grid.
 \end{itemize}

 After the initialization of the flow field is completed, ``\textit{Initialization}'' can re-output the flow field for inspection.

  \item The most important part of the software is the subroutine ``\textit{CalculateEigen}'', whose main idea is described in section \ref{subsection 2.2}. The role of this subroutine is to assemble the stability matrix and calculate its eigenvalues. To this end, it needs to call module ``\textit{Calculate Gradients}'' to calculate the gradients of the numerical flux, which can be divided into three parts as follows
  
  \begin{itemize}
	\item Firstly, as shown in (\ref{eq MUSCL method}), the variables on the left and right sides of the interface should be obtained by the module ``\textit{Reconstruct}''. As mentioned in section \ref{subsection 2.1}, there are two methods to reconstruct the variables: MUSCL and ROUND. And five different limiters including superbee \cite{Roe1985}, van Leer \cite{BramVanLeer1974}, van Albada \cite{VanAlbada1982}, minmod \cite{Roe1986a}, and the limiter proposed by Deng et al.\cite{Deng_Unified_2023} are contained in the MUSCL approach.
	
	\item The numerical flux can be obtained form the module ``\textit{Calculate Flux}'' by employing different Riemann solvers, including Roe \cite{Roe1981}, HLLC \cite{Toro1994}, HLL \cite{Harten1983}, van Leer \cite{VanLeer1997}, AUSM$ ^+ $ \cite{Liou1996}, SLAU \cite{Shima2011}, HLLE \cite{Einfeldt1988}, HLLEM \cite{Einfeldt1991}.
	
	\item Then, the gradients of the numerical flux can be calculated by the centered difference approximation.
  \end{itemize}

  \item Finally, ``\textit{Main Program}" will call the subroutine ``\textit{Output}'' to display the results, including the eigenvalues with the maximal real part, scatters of all the eigenvalues, and the unstable mode.
\end{itemize}
It should be noted that MSAT can run on both Windows and Linux operating systems. The versions of Windows and Linux have the same architecture.

\section{Illustrative examples}\label{section 3}

\subsection{Quantitative validation of MSAT}\label{subsection 3.1}

From equation (\ref{eq solution of linear error evolution model}), it can be found that the maximal real part of the eigenvalues can indicate the exponential growth rate of the initial perturbation error. As a result, the responsibility of MSAT can be quantitatively validated by comparing $ \max (\operatorname{Re}(\lambda(\mathbf{S})))  $ with the exponential growth rate, which is easy to obtain. The detail of this validation method can be found in \cite{Dumbser2004}. Note that the 2D flow field is initialized by the 1D computation throughout the paper if not mentioned specifically. 

\begin{figure}[htbp]
	\centering
	\includegraphics[width=0.6\textwidth]{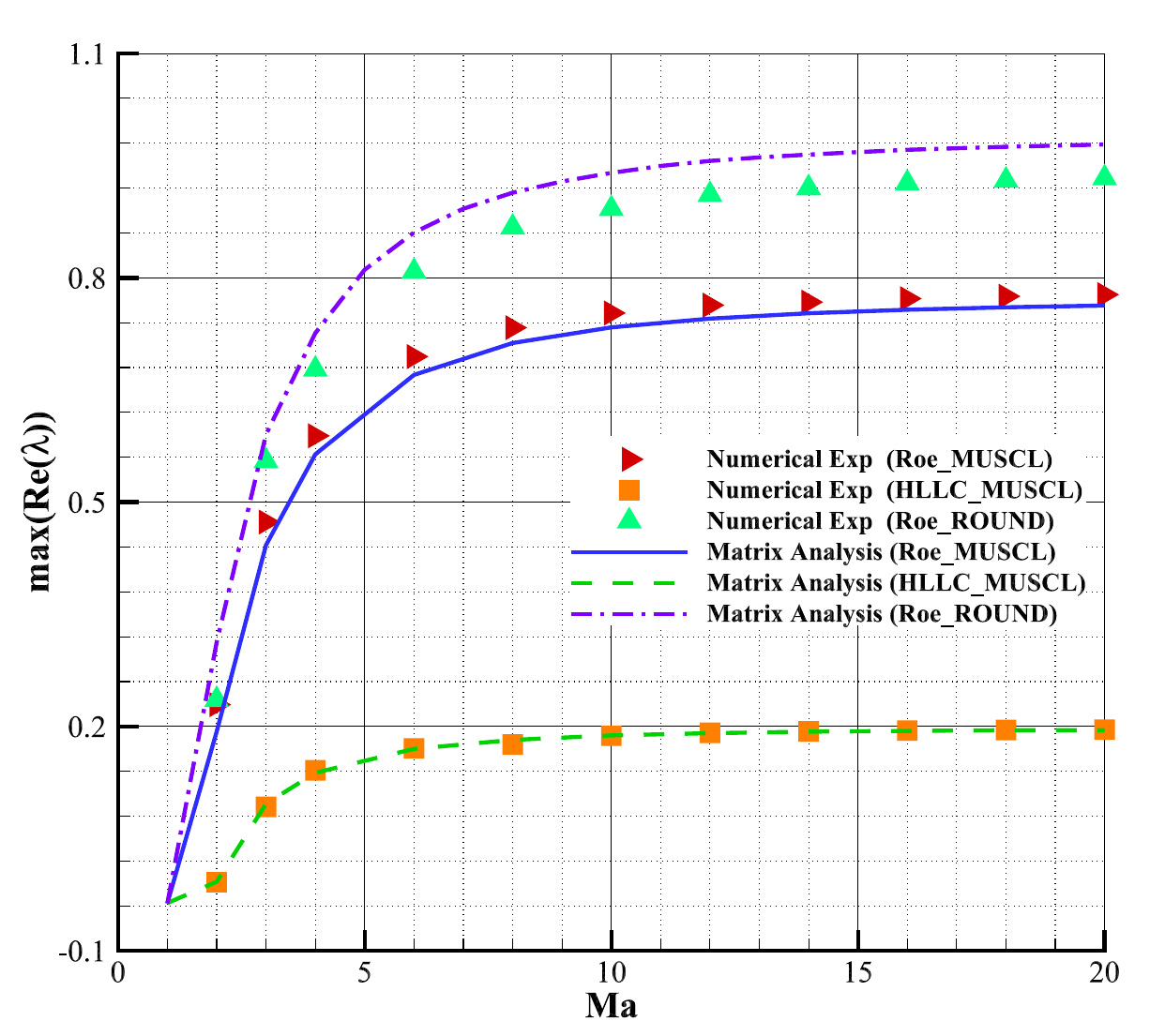}
	\caption{Quantitative validations of MSAT.(Grid with 11$ \times $11 cells, Roe and HLLC solvers, van Albada limiter in MUSCL, $ M_0=1,2,3,4,6,8,\cdots,20 $, and $ \varepsilon=0.1 $.)}
	\label{fig verify}
\end{figure}

Both the MUSCL approach (with the van Albada limiter) and ROUND scheme are employed to validate the respectively of MSAT. The Roe and HLLC solvers (with the wave speed estimates proposed by Davis \cite{Davis1988}) are used to compute the numerical flux. The conditions are $ M_0=1,2,3,4,6,8,\cdots,20 $, and $ \varepsilon=0.1 $. The computational grid is 11$ \times $11 Cartesian grid. The comparison between $ \max (\operatorname{Re}(\lambda(\mathbf{S})))$ and the exponential growth rate is shown in Fig.\ref{fig verify}. As shown, there are good agreements between them, and good agreements can also be obtained by other schemes with different solvers and limiter functions, which confirms the reliability of MSAT.
 
\begin{figure}[htbp]
	\centering
	\subfigure[First-order scheme with HLL solver.]{
	\begin{minipage}[t]{0.46\linewidth}
	\centering
	\includegraphics[width=0.9\textwidth]{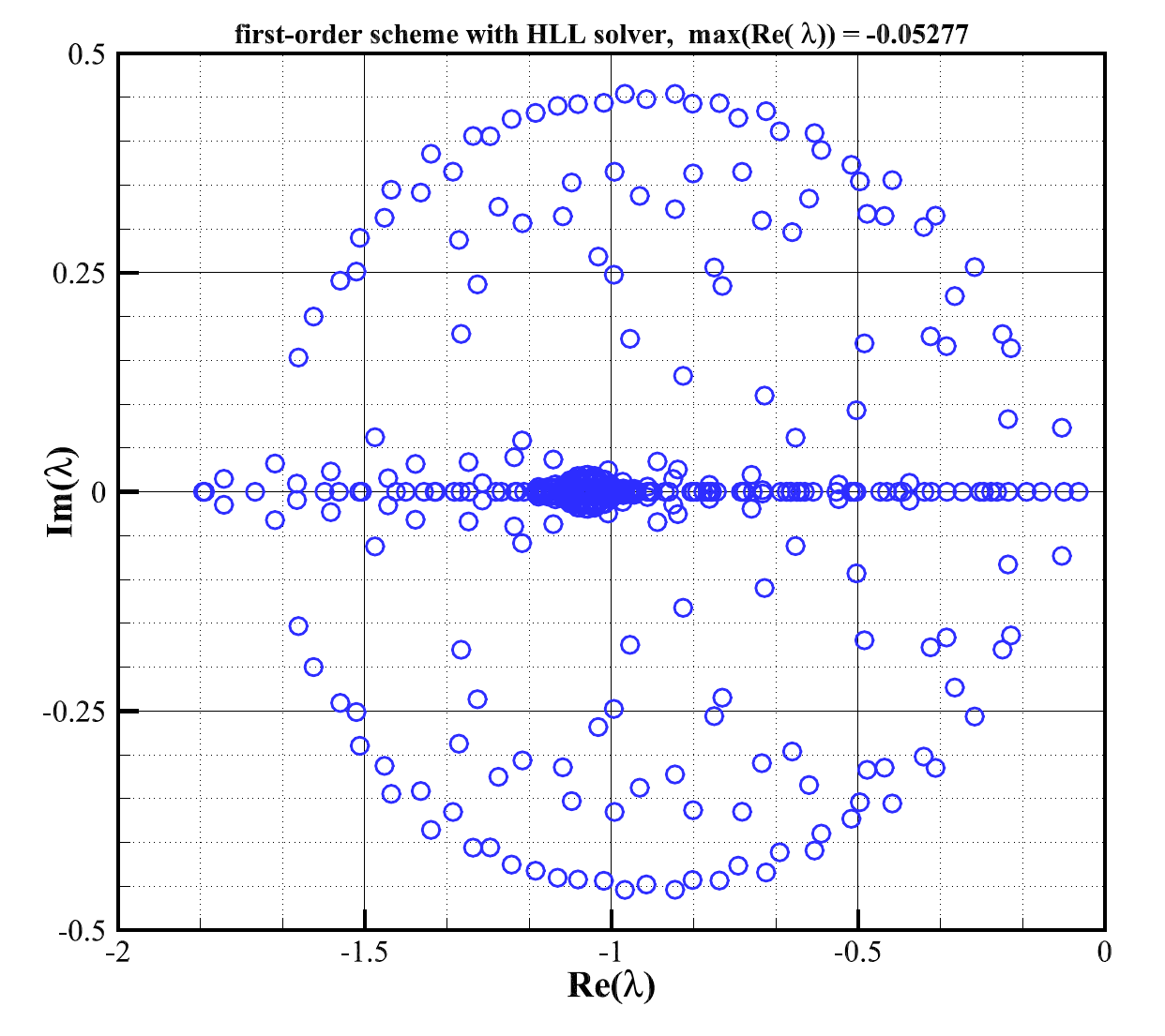}
	\end{minipage}
	}
	\subfigure[Second-order scheme with HLL solver.]{
	\begin{minipage}[t]{0.46\linewidth}
	\centering
	\includegraphics[width=0.9\textwidth]{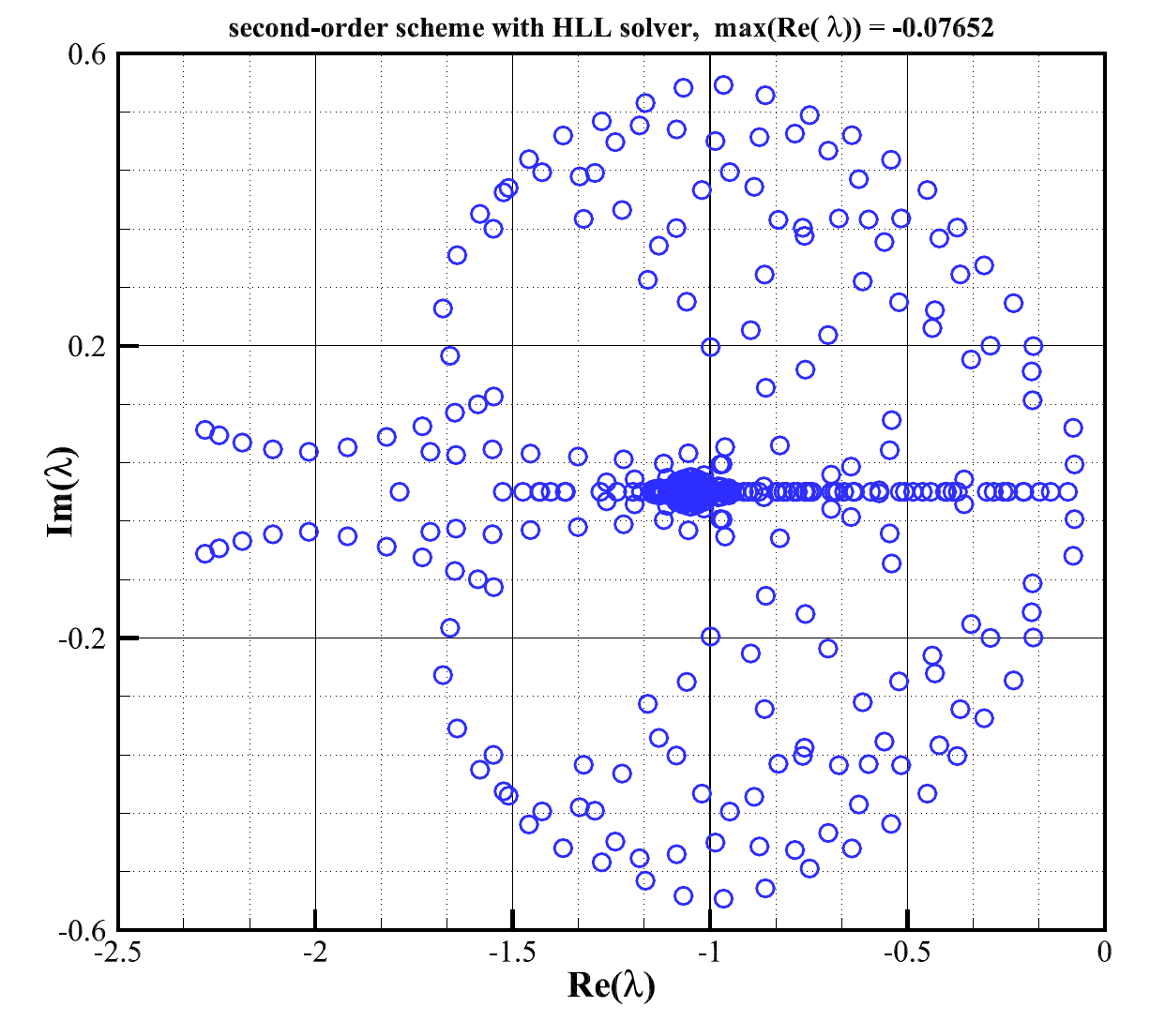}
	\end{minipage}
	}

	\subfigure[First-order scheme with HLLC solver.]{
	\begin{minipage}[t]{0.46\linewidth}
	\centering
	\includegraphics[width=0.9\textwidth]{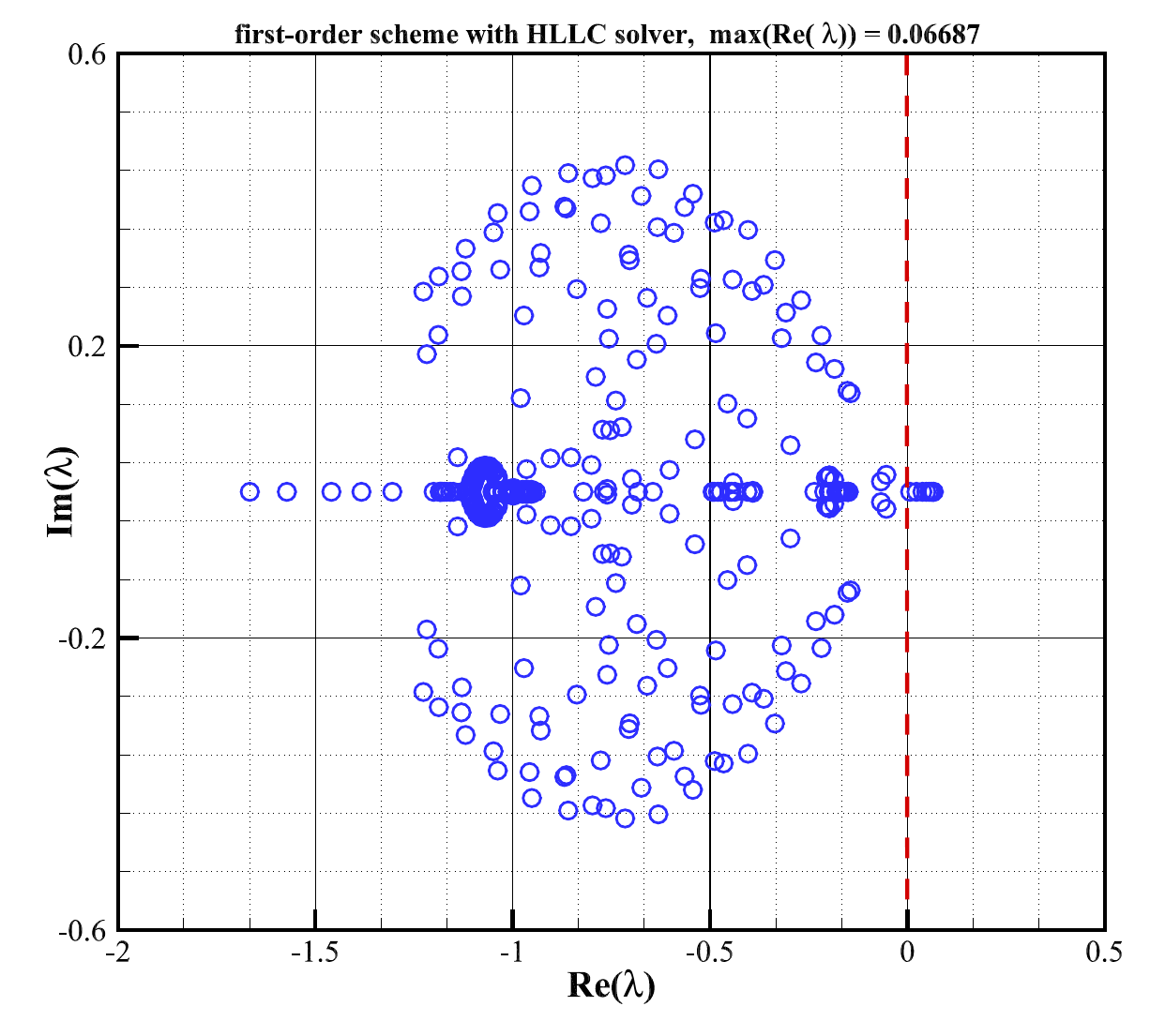}
	\end{minipage}
	}
	\subfigure[Second-order scheme with HLLC solver.]{
	\begin{minipage}[t]{0.46\linewidth}
	\centering
	\includegraphics[width=0.9\textwidth]{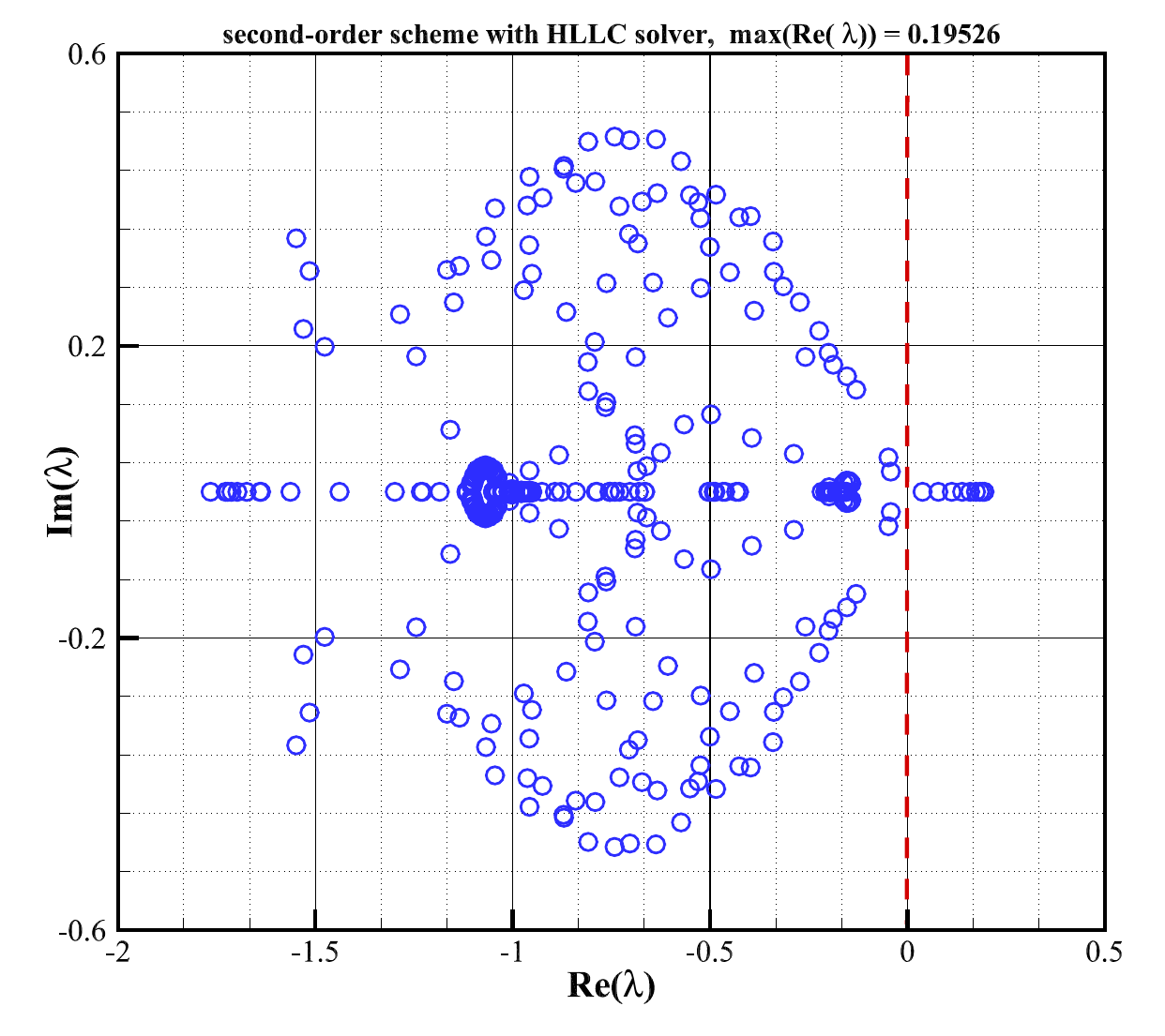}
	\end{minipage}
	}

	\caption{Distribution of the eigenvalues with different Riemann solvers and spatial accuracy.(Grid with 11$ \times $11 cells, van Albada limiter is used in the second-order scheme, $ M_0=20 $, and $ \varepsilon=0.1 $.)}\label{fig Riemann solver and spatial accuracy}
\end{figure}

\subsection{Influence of the Riemann solver and spatial accuracy on the shock instability}\label{subsection 3.2}

It has been demonstrated that the Riemann solver and spatial accuracy will significantly influence the shock instability\cite{Kitamura2009,Dumbser2004,Liou1996,Tu2014}. As shown in Fig.\ref{fig Riemann solver and spatial accuracy}, MSAT can be used to investigate how the two factors affect the shock instability problem. It can be found that the maximal real parts of all eigenvalues of the HLL solver are less than zero, indicating its stability. However, for the scheme using the HLLC solver, there are eigenvalues with positive real parts. Therefore, the computation will be unstable when simulating strong shocks with the HLLC solver. Moreover, as shown in Fig.\ref{fig Riemann solver and spatial accuracy}, the distribution of the eigenvalues for the two solvers are different between the first and second-order cases, indicating that the spatial accuracy can affect the shock instability for HLL and HLLC solvers. Specifically, the second-order scheme with the HLLC solver has a more significant maximal real part of the eigenvalues than the first-order case. In contrast, when employing the HLL solver, the second-order scheme has smaller real parts. So, it can be inferred that for the HLL solver, the stability of capturing strong shocks will be better when the spatial accuracy is enhanced to the second-order. However, the second-order scheme will be more unstable with the HLLC solver.

\begin{figure}[htbp]
	\centering
	\subfigure[The unstable mode for $ \rho $.]{
	\begin{minipage}[t]{0.46\linewidth}
	\centering
	\includegraphics[width=0.9\textwidth]{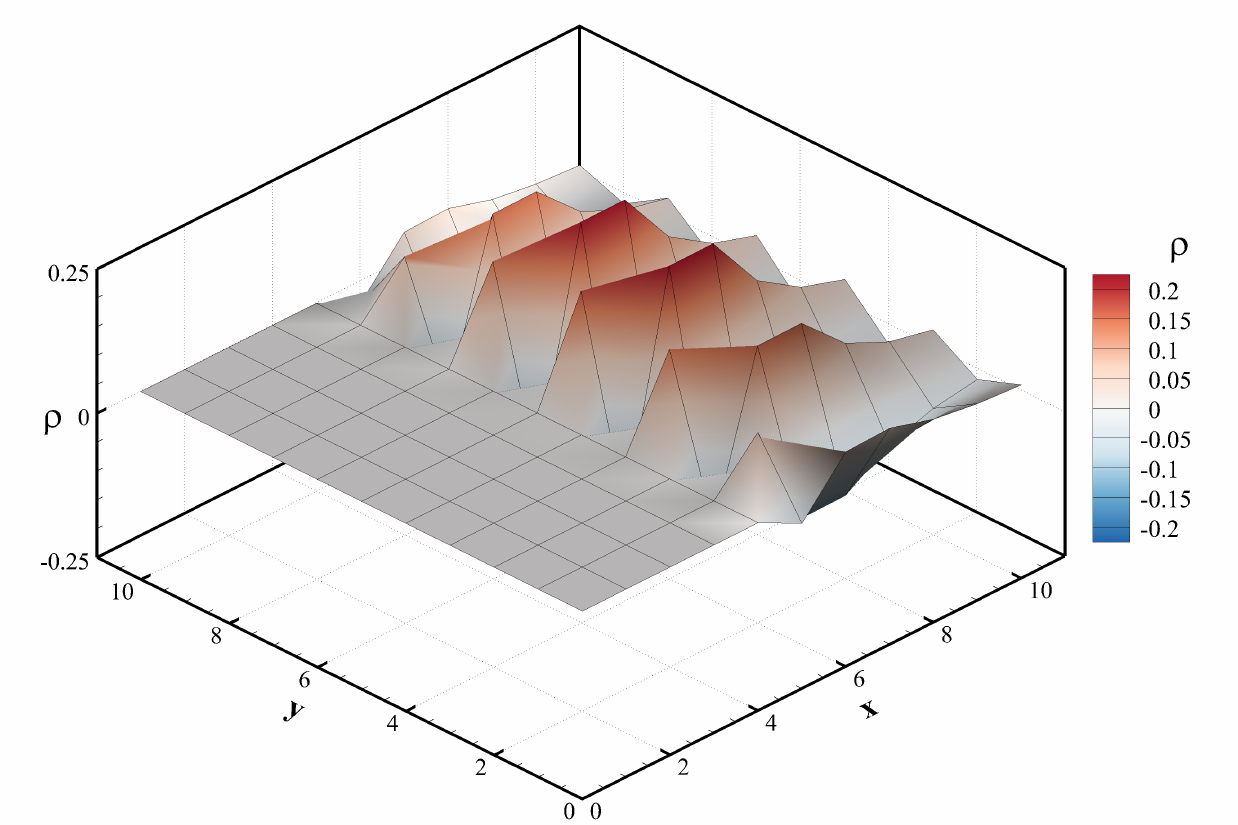}
	\end{minipage}
	}
	\subfigure[The unstable mode for u.]{
	\begin{minipage}[t]{0.46\linewidth}
	\centering
	\includegraphics[width=0.9\textwidth]{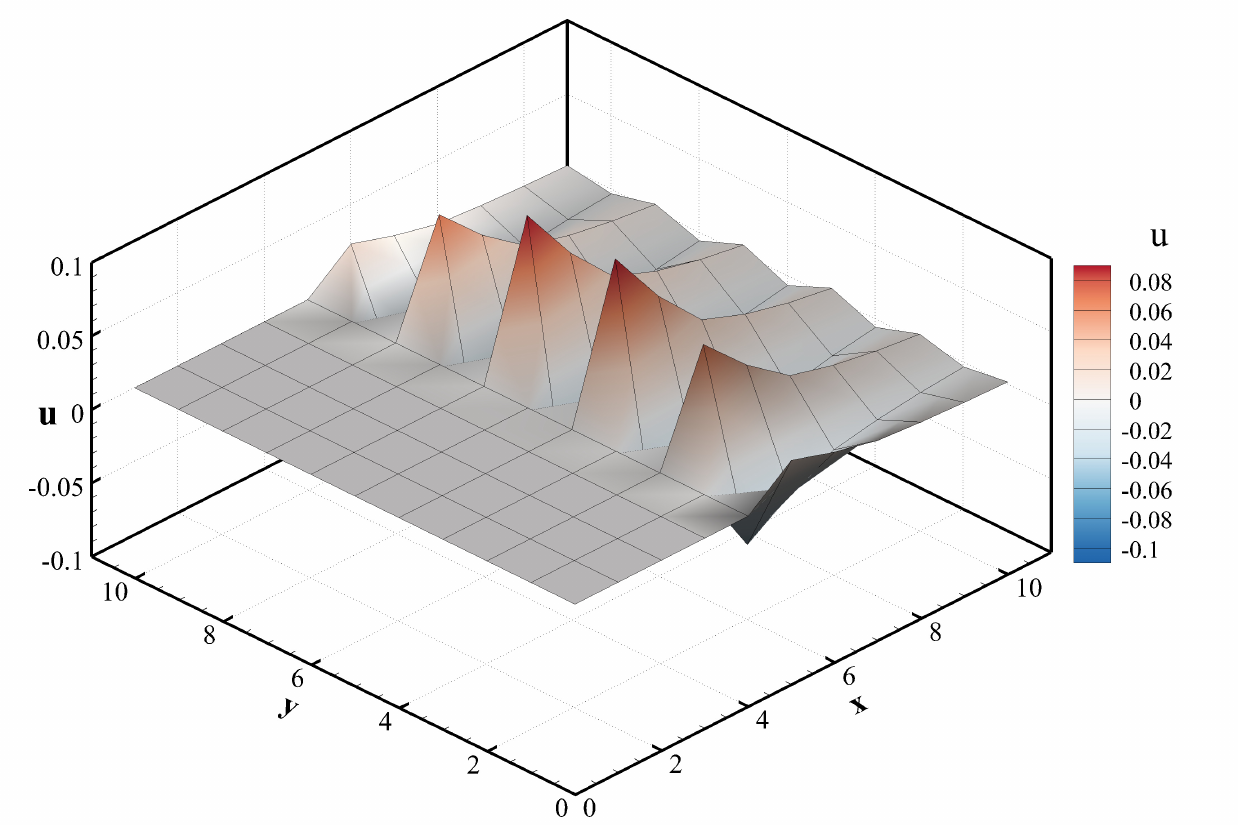}
	\end{minipage}
	}

	\subfigure[The unstable mode for v.]{
	\begin{minipage}[t]{0.46\linewidth}
	\centering
	\includegraphics[width=0.9\textwidth]{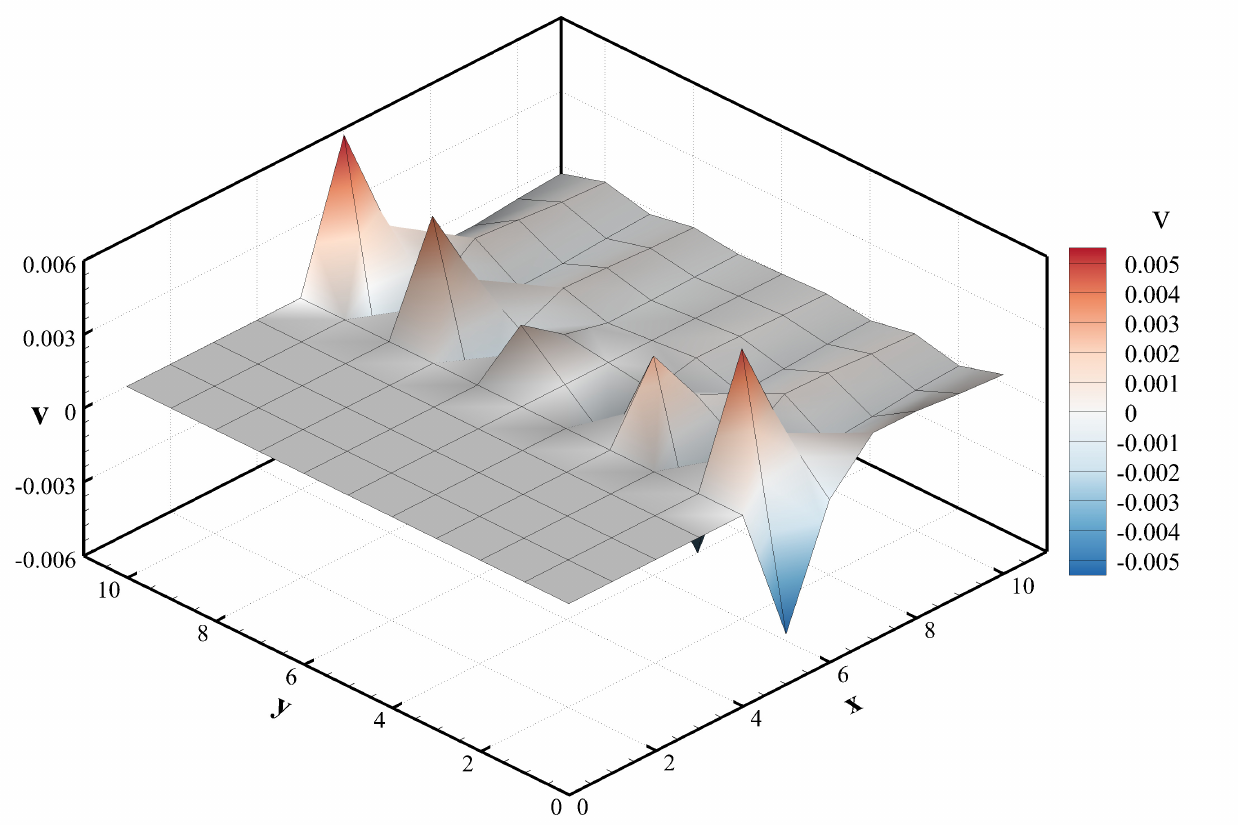}
	\end{minipage}
	}
	\subfigure[The unstable mode for p.]{
	\begin{minipage}[t]{0.46\linewidth}
	\centering
	\includegraphics[width=0.9\textwidth]{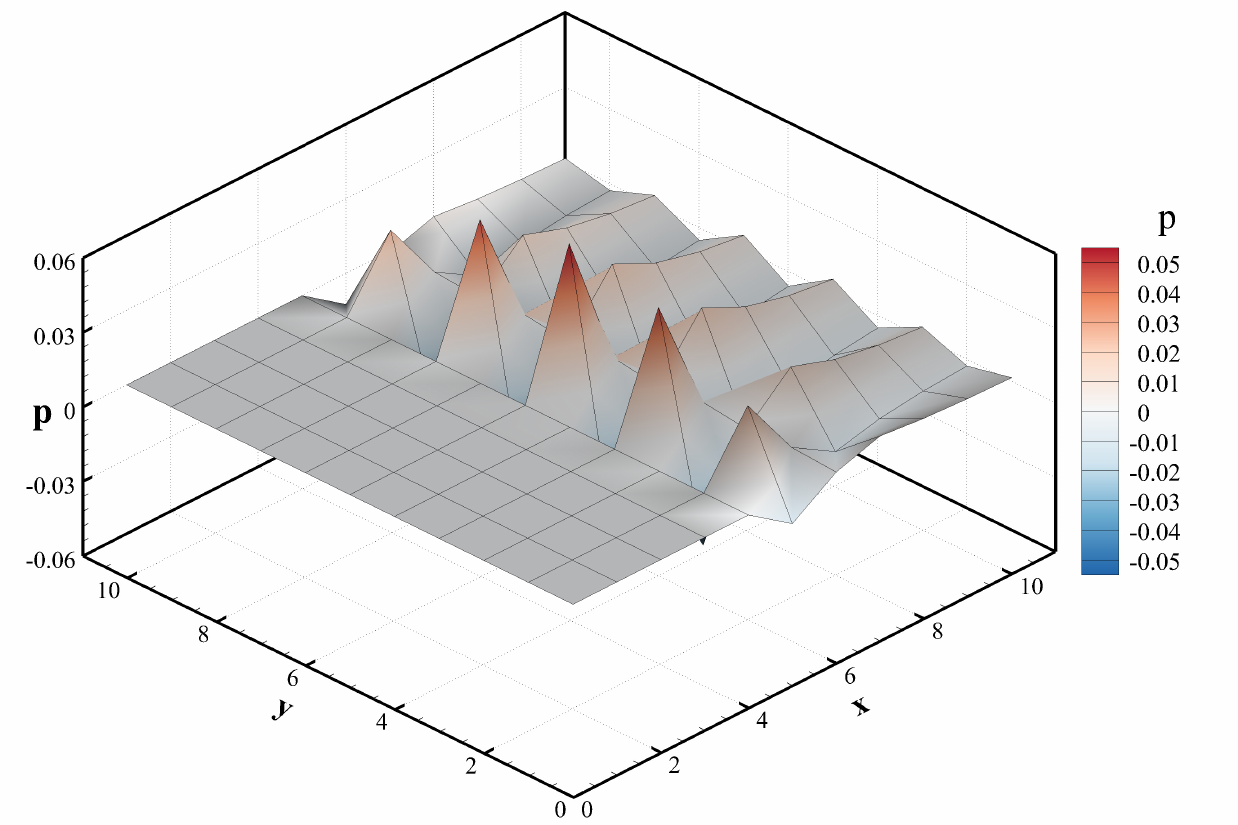}
	\end{minipage}
	}

	\caption{The unstable mode for $ \lambda = 0.19526+0i $.(Grid with 11$ \times $11 cells, second-order scheme with HLLC solver and van Albada limiter, $ M_0=20 $, and $ \varepsilon=0.1 $.)}\label{fig eigvector}
\end{figure}

This tool can also calculate the eigenvectors corresponding to the most unstable eigenvalue, which contains the information of the spatial behavior of the unstable mode.\cite{Dumbser2004} Fig.\ref{fig eigvector} shows the unstable mode of the second-order scheme with HLLC solver. As shown, the initial perturbation error is easier to influence the cells of $ j = 6 $, which contain the numerical shock structure. So, it can be inferred that the shock instability originates from the numerical shock structure, which is consistent with the conclusion in \cite{Xie2017}. Also, we can find that the downstream region is more susceptible to being affected by the unstable information from the numerical shock structure than the upstream region.

\begin{figure}[htbp]
	\centering
	\subfigure[Grid with 50$ \times $50 cells.]{
	\begin{minipage}[t]{0.46\linewidth}
	\centering
	\includegraphics[width=0.9\textwidth]{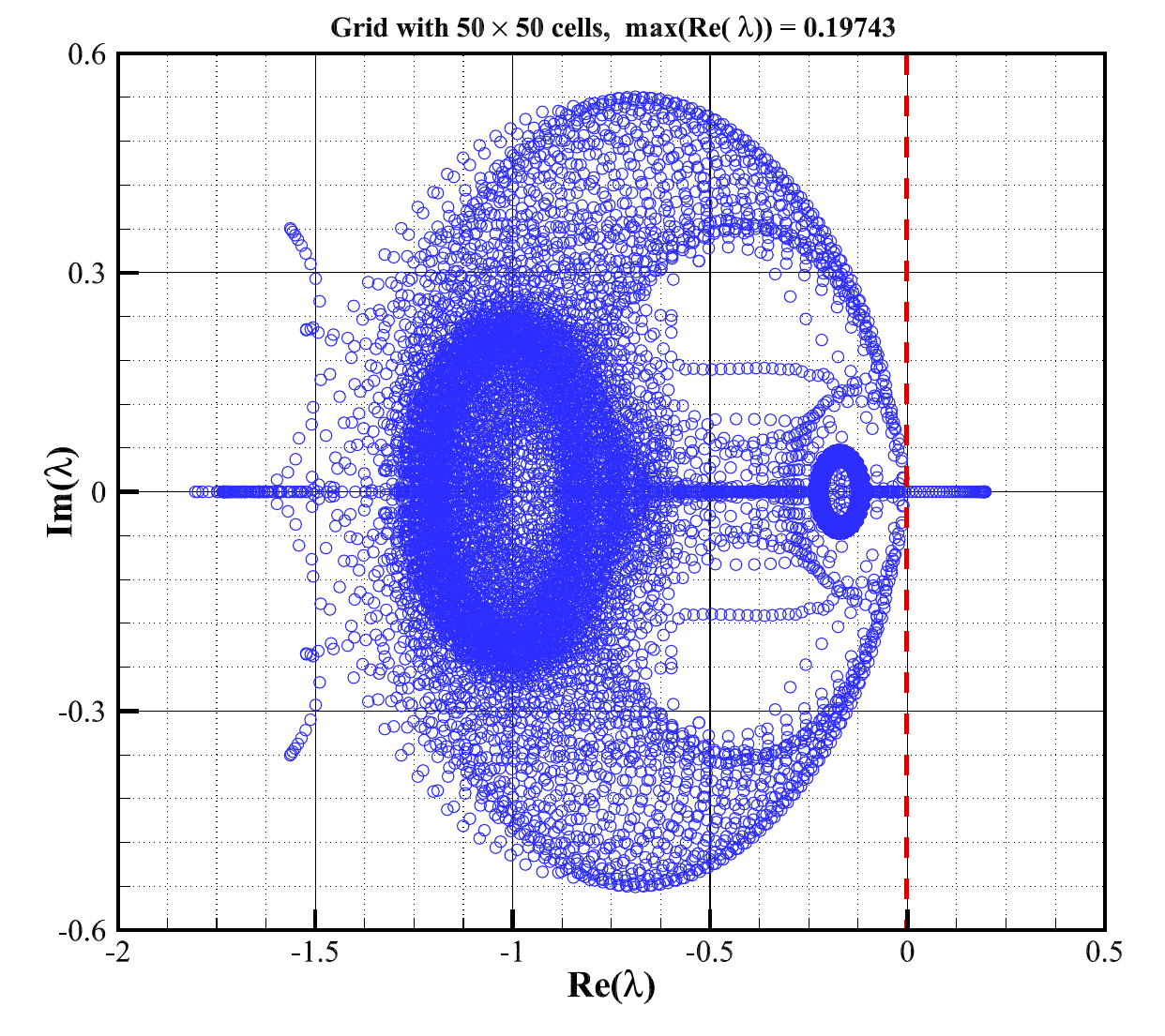}
	\end{minipage}
	}
	\subfigure[Grid with 50$ \times $10 cells.]{
	\begin{minipage}[t]{0.46\linewidth}
	\centering
	\includegraphics[width=0.9\textwidth]{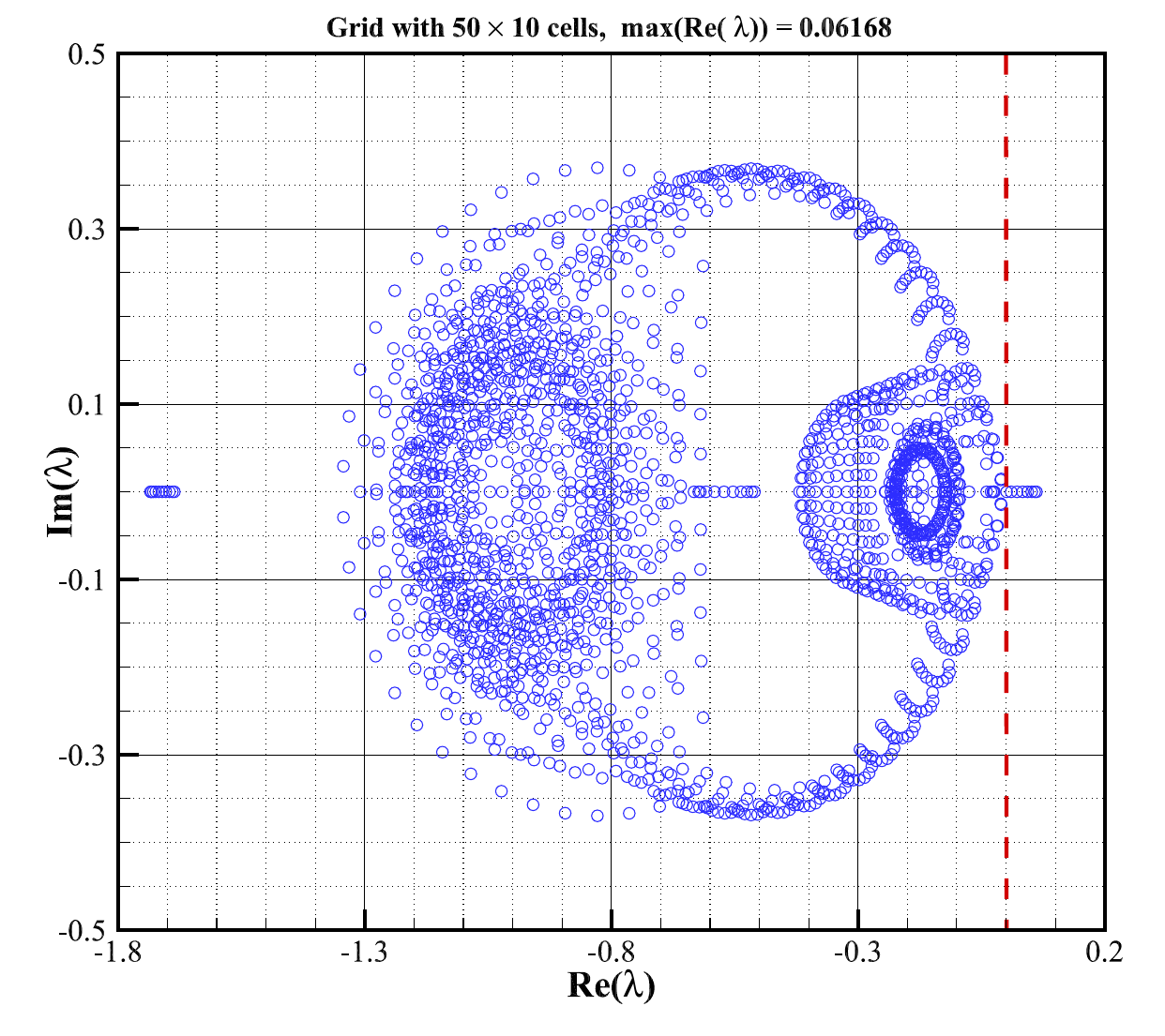}
	\end{minipage}
	}
	\caption{Influence of the computational grid with HLLC solver.(Grid with 50$ \times $50 and 50$ \times $10 cells respectively, second-order scheme with van Albada limiter, $ M_0=20 $, and $ \varepsilon=0.1 $.)}\label{fig Grid}
\end{figure}

\subsection{Influence of the computational grid on the shock instability} \label{subsection 3.3}

MSAT can also investigate the effect of the computational grid, which is also an important issue that will influence the shock instability \cite{Dumbser2004,Henderson2007,Ohwada2013}. Fig.\ref{fig Grid} shows all the eigenvalues computed with different computational grids and Fig.\ref{fig Grid}(a) employs the grid with 50$ \times $50 cells (grid A), while Fig.\ref{fig Grid}(b) is computed with 50$ \times $10 cells (grid B, 50 cells in x-direction and 10 cells in y-direction). Since the computational domain is 50$ \times $50 for the two grids, the aspect ratio $ \delta = \Delta y/ \Delta x$ differs, and $ \delta_A = 1, \delta_B =5$. As shown, as the aspect ratio increases, the maximal real part of the eigenvalues decreases, indicating that the shock instability can be alleviated by increasing the aspect ratio. This conclusion is consistent with the results of numerical experiments in \cite{Henderson2007,Ohwada2013}. Also, it can be found that the computation is still unstable even $ \delta_B =5 $, as the maximal real part of the eigenvalues still exceeds zero.

\subsection{Analysis of the hypersonic ﬂow over a cylinder}\label{subsection 3.4}
\begin{figure}[htbp]
	\centering
	\subfigure[HLL solver.]{
	\begin{minipage}[t]{0.46\linewidth}
	\centering
	\includegraphics[width=0.9\textwidth]{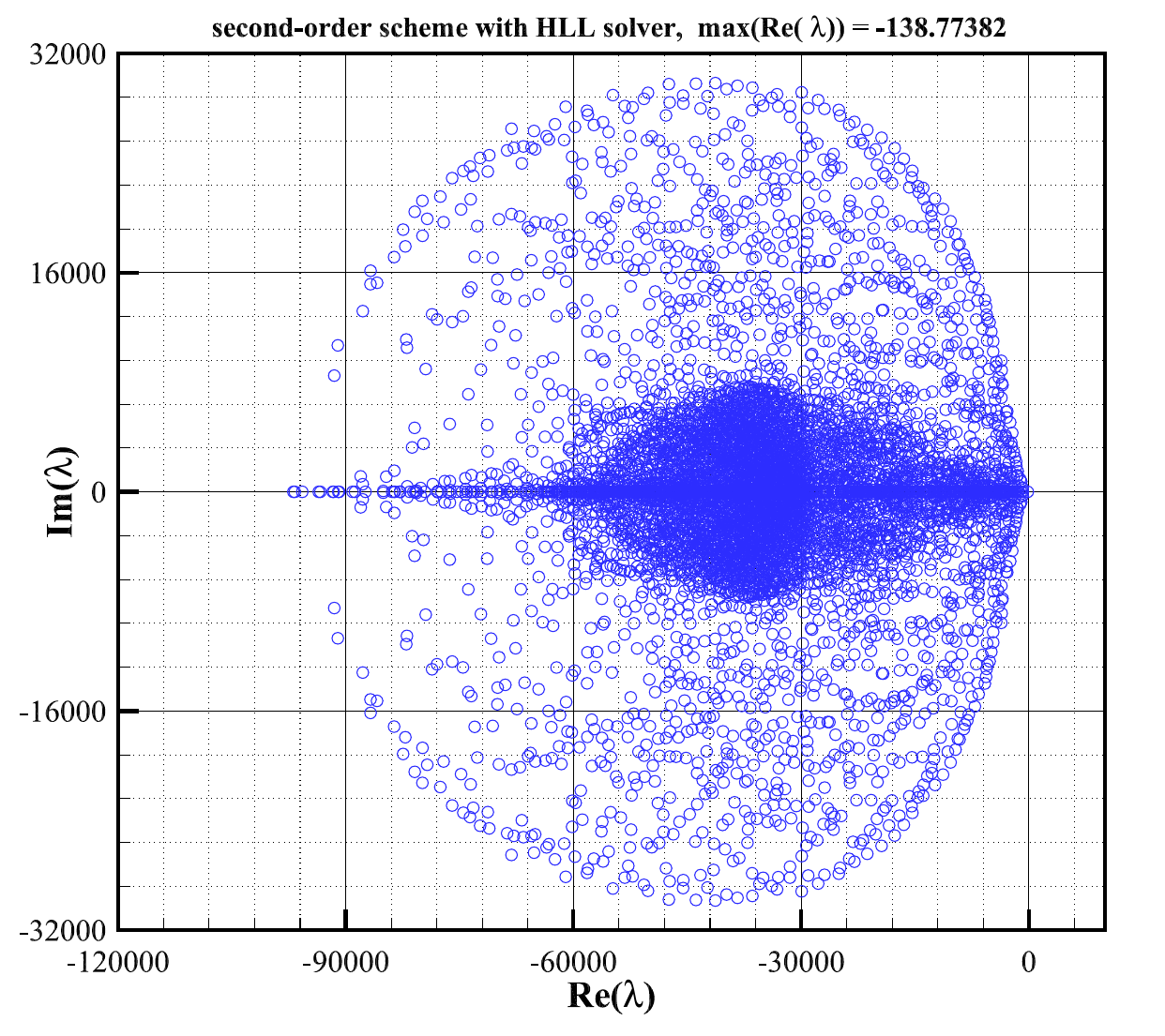}
	\end{minipage}
	}
	\subfigure[HLLC solver.]{
	\begin{minipage}[t]{0.46\linewidth}
	\centering
	\includegraphics[width=0.9\textwidth]{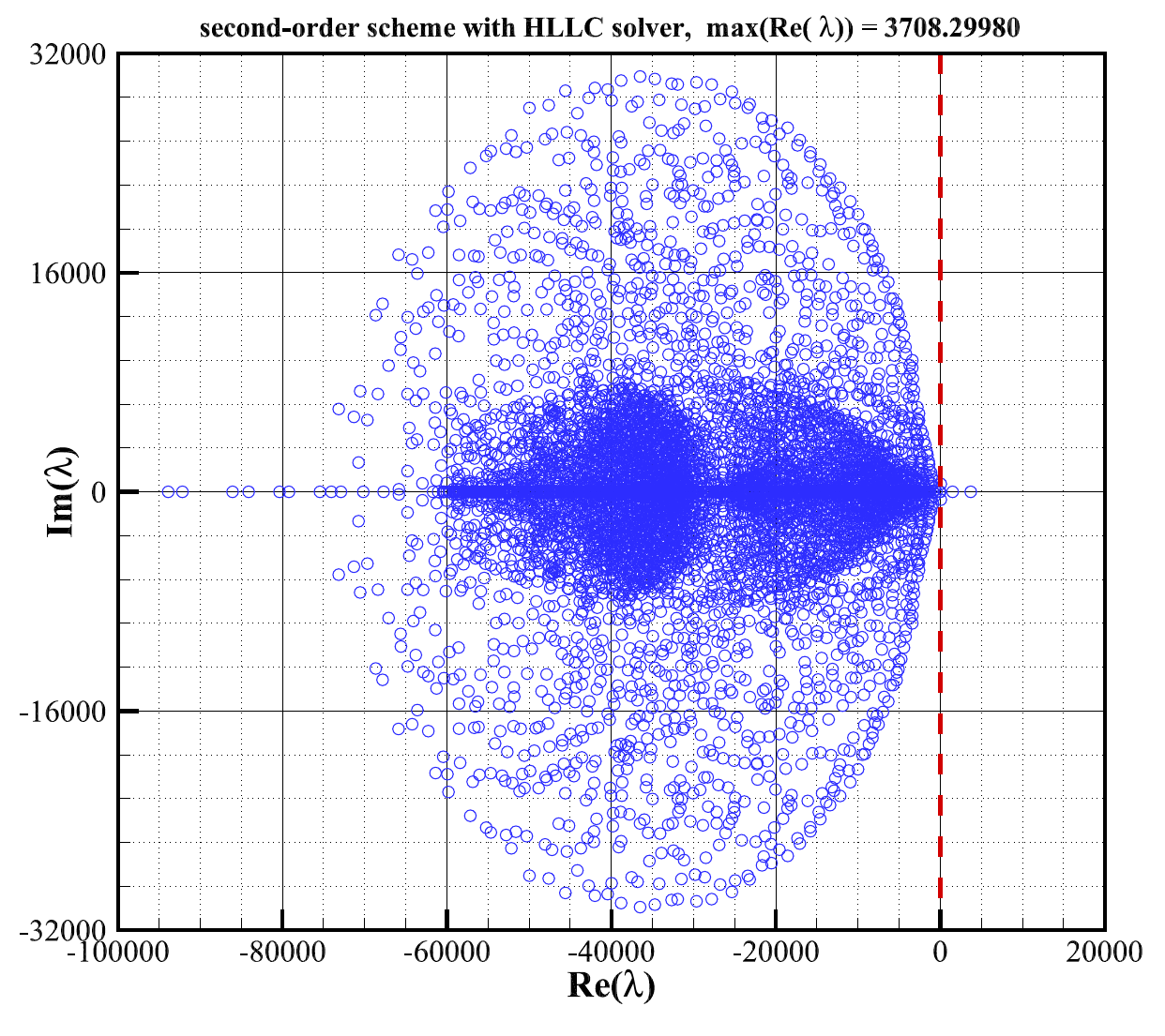}
	\end{minipage}
	}
	\caption{Matrix analysis of the hypersonic ﬂow over a cylinder.(Grid with 20$ \times $80 cells, second-order scheme with van Albada limiter.)}\label{fig cylinder}
\end{figure}

Illustrative examples of MSAT are displayed in section \ref{subsection 3.1} to \ref{subsection 3.3}. These examples are all based on the 2D steady normal shock problem, which is frequently employed in investigating the shock instability problem. Furthermore, MSAT has broader applicability and can be used to analyze the stability in practical calculations involving strong shocks. In this section, we analyze the hypersonic flow over a cylinder to illustrate its ability to handle other test cases. The computational grid is 20$\times$80 structured grid, and the conditions are as follows:
\begin{equation}   
	M_0=20.0, \quad \rho=1.4,  \quad p=1.0.
\end{equation} 
The initial flow was obtained by the second-order scheme with the HLL solver from the program solving two-dimensional Euler equations. Results are shown in Fig.\ref{fig cylinder}, from which it can be observed that all eigenvalues have negative real parts when the HLL solver is utilized. So computation of the hypersonic flow over a cylinder will be stable. However, if the HLLC solver is used, positive real parts can be observed, indicating its instability. As noted in \cite{Simon2018a,Pandolfi2001}, the HLLC solver will suffer from the shock instability problem while simulating the hypersonic flow over a cylinder, whereas the HLL solver can yield stable results. The results from the matrix analysis in this study are consistent with the conclusion obtained by \cite{Simon2018a,Pandolfi2001}, thereby confirming the accuracy of MSAT in analyzing practical test cases.

\section{Impact and conclusion}\label{section 4}
The shock instability problem severely limits the application of shock-capturing methods in supersonic or hypersonic flow simulations. Currently, there is a lack of open-source tools for quantitative analyzing the shock instability problem of the higher-order scheme. Based on the work in \cite{2305.03281}, we develop the unified matrix stability method for the schemes with three-point stencils and present an open-source tool, MSAT, in this paper. MSAT is a helpful tool to quantitatively analyze the shock instability problem. The responsibility of MSAT has been well verified by numerical experiments. Two reconstruction methods, MUSCL and ROUND, and widely used Riemann solvers and are integrated into MSAT. Also, it is simple to add new reconstruction methods and solvers in the code if needed.

Illustrative examples in section \ref{section 3} demonstrate some applications of MSAT. In addition, MSAT has a broader spectrum of applications. Apart from analyzing the impact of spatial accuracy, Riemann solvers, and computational grid on the shock instability problem, MSAT can be directly used to analyze more issues, including limiter functions, reconstruction methods, shock intensity, and numerical shock structure. Meanwhile, the underlying mechanism of shock instabilities can also be investigated by this tool. Readers are referred to references \cite{Dumbser2004,Shen2014} for the similar study. Furthermore, MSAT also has the ability to evaluate whether the practical simulation of supersonic/hypersonic flows will suffer from the shock instability problem. Moreover, with the comprehensive understanding of the framework of MSAT and necessary modifications, MSAT can be integrated as a module into the program solving Euler equations, NS equations, shallow water equations, and so on. Before the computation starts, a simple analysis can be performed by MSAT to guide the choice of numerical schemes. In this way, the shock instability problem can be avoided.

As a result, with the help of MSAT, there is no need for researchers to build the matrix stability analysis program from scratch. They can devote more time and effort to investigating the shock instability phenomenon, contributing to a better understanding of the shock instability problem and shedding new light on developing reliable shock-capturing schemes. MSAT is part of our series of works. Currently, we are conducting research on the matrix stability analysis method for high-order schemes, which will be integrated into MSAT in the future. In this way, MSAT can also be used to analyze the shock instability problem for high-order schemes.

\section*{Declaration of competing interest}\label{section 5}
The authors declare that they have no known competing financial interests or personal relationships that could have appeared to influence the work reported in this paper.

\section*{Acknowledgements}	
The authors appreciate reviewers’ useful comments and valuable suggestions for the original manuscript. This work was supported by National Natural Science Foundation of China (Grant No.12202490), Natural Science Foundation of Hunan Province, China (Grant No. 11472004), the Scientific Research Foundation of NUDT (Grant No. ZK21-10), and Postgraduate Scientific Research Innovation Project of Hunan Province (Grant Nos. CX20220010 and CX20220036).

\bibliography{mybibfile}

\end{document}